\documentclass{article}
 
\usepackage{amssymb,amsmath,amsthm} 
\usepackage{graphicx}
\usepackage[latin1]{inputenc}
\usepackage[a4paper, centering, text={6.5in, 10in}]{geometry} 
\usepackage{authblk}
\usepackage{cancel}

\newtheorem{theorem}{Theorem}[section]

\newtheorem{proposition}[theorem]{Proposition}
\newtheorem{corollary}[theorem]{Corollary}
\newtheorem{definition}[theorem]{Definition}
\newtheorem{example}{Example}[section]
\newtheorem{remark}[example]{Remark}

\newcommand{\s}{\sigma}

\providecommand{\keywords}[1]{\textbf{\textit{Index terms---}} #1}

\begin{titlepage}

\author{Francisco Mota}
\affil{Departamento de Engenharia de Computação e Automação\\
Universidade Federal do Rio Grande do Norte -- Brasil\\
e-mail:mota@dca.ufrn.br}

\date{\today}

\title{A Note on Gamma Function and Convolution}

\end{titlepage}

\begin{document}

\maketitle

\begin{abstract}
In this note we explore the relationship between the operation of convolution of functions 
and the Eulerian integrals. This approach allow us to obtain some expressions for the convolution of a 
certain class of functions in terms of the Gamma Function as well as to derive some well known properties 
of the Gamma Function  by using the concept and properties of the convolution.


\end{abstract}

\section{Introduction}

Gamma function dismisses any introduction due to its importance in all mathematics; its properties have been 
analysed in literature for a long time. This note is, in a certain sense, inspired by the classical 
monograph \cite{artin} in which the author obtain some important properties  of Gamma function by 
using an argument of ``unicity''. We follow his approach to show the unicity of a certain 
set of functions that satisfies some properties (similar to the ones proposed in \cite[pp.~13]{artin} for Gamma function);
in the following we use the relationship of convolution with the Eulerian integrals to obtain some identities 
involving convolution of functions in this set that resembles the classical ones for the 
Gamma function. In particular, we obtain a ``convolutional'' version of the so-called Gauss product formula 
in Proposition~\ref{convgprodt}, which, in turn, implies the classical product formula. Additionally, 
in Propostion~\ref{convinv}, we use the fact that the domain of Gamma function may be extended to all the set of Reals 
in order to obtain a formula for the convolution of signals even when the convolution integral diverges.

\section{Convolution of Functions}\label{adconv}

Let be a function which maps the set of Reals ($\mathbb{R}$) into the set of of Complex ($\mathbb{C}$),
and it is defined almost everewhere on its domain:
\begin{align*}
f \colon & \mathbb{R} \to \mathbb{C}\\
  & t \mapsto f(t)
\end{align*}
Among the operations we can realize with functions, one of the most important is {\em convolution}~\cite{wiki}, which 
is a (binary) operation as defined below:

\begin{definition}\label{aconv}\em 
Consider two functions $f$ and $g$, its convolution, representd by ``$f*g$", is defined as:
\begin{equation}\label{aconveq}
(f*g)(t) = \int_{-\infty}^{\infty}f(\tau)g(t-\tau)d\tau
\end{equation} 
\end{definition}

\begin{remark}\em 
In order the convolution to be well defined as an operation we require that integral in Equation~(\ref{aconv}) to be 
absolutely convergent, and under this condition, convolution is  a commutative 
and associative binary operation. On the other hand, there is an issue related to the existence 
of the identity for this operation, and, in fact, it is well known that the identity for convolution is 
not a usual function, but a distribution \cite{wiki}.

\end{remark}

\begin{remark}\em If both functions , $f$ and $g$, are such that $f(t)=g(t)=0$ for $t<0$, the convolution integral 
(\ref{aconveq}) becomes null for $t<0$, and for $t>0$ we get:
\begin{equation}\label{aconveq0}
(f*g)(t) = \int_{0}^{t}f(\tau)g(t-\tau)d\tau
\end{equation} 

\end{remark}

\begin{remark}\em The $n$-fold convolution of a function with itself is generally respresented by $f^{*n}$, 
that is
\[
f^{*n} = \underbrace{f*f*\cdots *f}_{n\text{ times}}
\]
\end{remark}

\subsection{Convolution Derivative Property, Units Step Function and its Derivative}\label{basics}

One important property of convolution is that its derivative depends only on the derivative of one of the functions 
involved in the operation \cite{wiki}; in fact, if $f$ and $g$ are functions, with $f$ or $g$ differentiable 
(i.e. $\dot f$ and/or $\dot g$ exists), then
\begin{equation}
\dot{(f*g)} = \dot f*g = f*\dot g \label{convdif}
\end{equation}
Additionaly, this property may be used to give an interpretation of the unity of the convolution operation as the 
derivative of the ``unit step function" as discussed in the following. 

The unit step function (represented by ``$\s$'') is defined as

\begin{equation}\label{ustep}
\s(t) = \begin{cases} 0, & t<0 \\ 1 & t > 0 \end{cases}
\end{equation}
If we consider the convolution of $\s$ with a function $g$ we obtain:
\begin{eqnarray}
(\s*g)(t) & = & \int_{-\infty}^{\infty}g(\tau)\s(t-\tau)d\tau, \text{ since}
\quad \s(t-\tau) = \begin{cases} 0, & \tau>t \\ 1 & \tau\le t \end{cases} \quad\text{we get,}\nonumber \\
(\s*g)(t) & = & \int_{-\infty}^{t}g(\tau)d\tau. \label{sgint}
\end{eqnarray}
Since the convolution with a unit step has the effect of integrating 
the other one function, as shown in (\ref{sgint}), if we consider a set of (test) differentiable functions 
$g$, and using the derivative property of convolution (\ref{convdif}), to ``transfer" the derivative 
from $\s$ to $g$, we may write:
\[
(\dot \s*g)(t) \triangleq (\s*\dot{g})(t) \underset{\text{by (\ref{sgint})}}{=} 
\int_{-\infty}^{t}\frac{d}{d\tau}g(\tau)d\tau = g(t),
\]
that is
\[
\dot{\s}*g \triangleq \s*\dot g = g.
\]
Then  ``$\dot \s$'' may have an interpretation  as the unity element for the operation 
of convolution, which is well-known to be the ``unit impulse" (``$\delta$"). So, we may write
\[
\dot \s =  \delta
\]

\section{The Set of Functions \boldmath $f_{\alpha}$}

In this note we consider a set of of real valued functions  parameterized by a real positive parameter. 
In order to characterize this set of functions, we follow an approach presented in \cite[pp.~13]{artin} for 
the Gamma function, as shown below:

\begin{proposition}\em\label{talfa}
Let it be a set of real and non-negative (parameter dependent) functions ``$f_{\alpha}$", $\alpha>0\in\mathbb{R}$, 
\begin{align*}
f_{\alpha} \colon & \mathbb{R} \to \mathbb{R^+}\\
  & t \mapsto f_{\alpha}(t)
\end{align*}
which satisfies the following properties:
\begin{enumerate}
\item[(1)] For all $\alpha>0$, $f_{\alpha}(t)=0$ for $t<0$ and $f_{\alpha}(t)>0$ for $t>0$. 
\item[(2)] $f_1=\s$; that is $f_1$ is the unit step function, as defined in (\ref{ustep});
\item[(3)] $f_{\alpha+1} = \alpha f_{\alpha}*f_1 \iff \displaystyle f_{\alpha+1}(t) = \alpha\int_0^t f_{\alpha}(\tau)d\tau$
$\iff \dot{f}_{\alpha+1}(t) = \alpha f_{\alpha}(t)$.
\item[(4)] $f_{\alpha}(t)$, $t>0$, is log-convex \cite{artin} relative to parameter $\alpha$, that is,
\[
\ln f_{[\lambda\alpha_1+(1-\lambda)\alpha_2]}(t) \le \lambda \ln f_{\alpha_1}(t) + (1-\lambda)\ln f_{\alpha_2}(t), 
\quad 0\le\lambda\le 1.
\]
\end{enumerate}

Then this set of functions is necessarily defined as:

\begin{equation}\label{fa}
f_{\alpha}(t) = \begin{cases} 0, & t<0 \\ t^{\alpha-1}, &t>0 \end{cases}
\end{equation}

\begin{proof} Of course, we only need to prove for $t>0$. Before all, we note that property (3) above, immediately 
implies the following results:
\begin{equation}\label{fni}
f_n(t)  = t^{n-1}, \quad t>0, \; n\in\mathbb{N}
\end{equation}
and
\begin{equation}\label{fan}
f_{\alpha+n}(t) = [(\alpha+n-1)(\alpha+n-2)\cdots(\alpha+1)\alpha](f_{\alpha}*f_1^{*n})(t), \quad t>0,
\end{equation}
where  $f_1^{*n}$ is the $n$-fold convolution ``$f_1*f_1*\cdots*f_1$'', and $f_1=\s$ is the unit step function.
Additionally we can differentiate Equation~(\ref{fan}) $``n"$ times and, by using the fact that $\dot{f}_1=\delta$,
along with derivative property (\ref{convdif}) of convolution, to obtain
\begin{equation}\label{dnfan}
f_{\alpha+n}^{(n)}(t) = [(\alpha+n-1)(\alpha+n-2)\cdots(\alpha+1)\alpha]f_{\alpha}(t), \quad t>0.
\end{equation}

Now we follow a similar approach presented in \cite[pp.~13]{artin}: By the fact that $f_{\alpha}(t)$, $t>0$, 
is log-convex relative to $\alpha$, we have that differences quocient in ``$\alpha$'' 
increase monotonically for $f_{\alpha}(t)$. Then by considering $0<\alpha\le 1$ and $n\ge 2$ we may write:
\begin{align}
\frac{\ln f_{-1+n}(t)-\ln f_n(t)}{(-1+n)-n}  &\le \frac{\ln f_{\alpha+n}(t)-\ln f_n(t)}{(\alpha+n)-n} 
\le \frac{\ln f_{1+n}(t)-\ln f_{n}(t)}{(1+n)-n} \implies\nonumber \\
\frac{\ln(t^{-1+n-1})-\ln(t^{n-1})}{(-1+n)-n} &\le \frac{\ln f_{\alpha+n}(t)-\ln f_n(t)}{(\alpha+n)-n} 
\le \frac{\ln(t^n)-\ln(t^{n-1})}{(1+n)-n} \implies\nonumber \\
(n-1)\ln t - (n-2)\ln t & \le \frac{\ln f_{\alpha+n}(t)-\ln f_n(t)}{\alpha}
\le n\ln t-(n-1)\ln t \implies\nonumber \\
\ln t & \le \frac{\ln f_{\alpha+n}(t)-\ln f_n(t)}{\alpha} 
\le \ln t \implies 
\frac{\ln f_{\alpha+n}(t)-\ln f_n(t)}{\alpha}  =  \ln t \label{lnf}
\end{align}
And so, develloping (\ref{lnf}) further we have 
\begin{align}
\ln f_{\alpha+n}(t) & = \alpha\ln t + \ln f_n \implies\nonumber \\
\ln f_{\alpha+n}(t) & = \ln t^{\alpha} + \ln f_n \implies\nonumber \\
\ln f_{\alpha+n}(t) & = \ln (t^{\alpha}f_n), \quad f_n(t) = t^{n-1} \implies\nonumber \\
\ln f_{\alpha+n}(t) & = \ln(t^{\alpha+n-1}) \implies\nonumber \\
f_{\alpha+n}(t)     & = t^{\alpha+n-1}, \quad t>0, \quad n\ge 2 \label{fan2}
\end{align}
We can differentiate (\ref{fan2}) above $``n"$ times to obtain
\begin{equation}\label{dnfan2}
f_{\alpha+n}^{(n)}(t) = [(\alpha+n-1)(\alpha+n-2)\cdots(\alpha+1)\alpha]t^{\alpha-1}.
\end{equation}
Comparing (\ref{dnfan2}) with (\ref{dnfan}) we get the result (\ref{fa}). 
\end{proof}

\end{proposition}

\begin{corollary}\em\label{talfac} Let be the $n$ differentiable set of functions $g_\alpha(t):\mathbb{R}\to\mathbb{R^+}$, 
such that $g_1$ is defined as
\begin{equation}\label{g1}
g_1=af_1^{*(n+1)}, \quad f_1^{*(n+1)}=\underbrace{f_1*f_1*\cdots*f_1}_{n+1 \text{ times}}
\end{equation}
where, $a>0\in\mathbb{R}$, $n\in\mathbb{N}$, $f_1=\sigma$ is the unit step function. 
If $g_{\alpha}$ satifies conditions (1), (3) and (4) of Proposition~\ref{talfa}, 
listed  for functions $f_{\alpha}$, Then:
\begin{equation}\label{ga}
g_{\alpha} = af_1^{*n}*f_{\alpha}.
\end{equation}

\begin{proof}
Since the $\dot{f}_1=\delta$, the $n$-derivativative of $f_1^{*(n+1)}=f_1$, we have
\[
g_1^{(n)}/a = f_1
\]
So ``$g_{\alpha}^{(n)}/a$'' satifies all four conditions (1), (2), (3) and (4) listed above for $f_{\alpha}$, and
then $g_{\alpha}^{(n)}/a=f_{\alpha}$ by Proposition~\ref{talfa}; that is
\begin{align*}
g_{\alpha}^{(n)}/a          & = f_{\alpha} \implies \\
f_1^{*n}*g_{\alpha}^{(n)} & = af_1^{*n}*f_{\alpha} \implies \\[0.2cm]
               g_{\alpha} & = af_1^{*n}*f_{\alpha}                 
\end{align*}

\end{proof}

\end{corollary}



\section{The Eulerian Integrals and Convolution}

The First Eulerian Integral  \cite{artin}, also known as the (real) Beta function $\mathcal{B}$, is defined as
\begin{equation}\label{fei}
\mathcal{B}(\alpha_1,\alpha_2) = \int_0^1x^{\alpha_1-1}(1-x)^{\alpha_2-1}dx, \quad \alpha_1>0,\alpha_2>0
\end{equation}

\begin{proposition}\em \label{fab}
Let be $f_{\alpha}$ as defined in (\ref{fa}), then:
\begin{equation}\label{fabe}
(f_{\alpha_1}*f_{\alpha_2})(t) = \mathcal{B}(\alpha_1,\alpha_2)f_{\alpha_1+\alpha_2}(t)
\end{equation}

\begin{proof} By using the definition of convolution in (\ref{aconveq0}), we have:
\begin{eqnarray*}
(f_{\alpha_1}*f_{\alpha_2})(t) & = & \int_0^t\tau^{\alpha_1-1}(t-\tau)^{\alpha_2-1}d\tau, \quad \tau=tx,\; t>0\\
                               & = & \int_0^1(tx)^{\alpha_1-1}(t-tx)^{\alpha_2-1}tdx \\
                               & = & t^{\alpha_1+\alpha_2-1}\int_0^1x^{\alpha_1-1}(1-x)^{\alpha_2-1}dx\\
                               & = & f_{\alpha_1+\alpha_2}(t)\mathcal{B}(\alpha_1,\alpha_2)
\end{eqnarray*}

\end{proof}

\end{proposition}

The Second Eulerian Integral \cite{artin}, also known as the (real) Gamma function $\Gamma$, is defined by 
\begin{equation}\label{sei}
\Gamma(\alpha) = \int_0^{\infty}t^{\alpha-1}e^{-t}dt, \quad \alpha>0.
\end{equation}

\begin{proposition}\label{gammaconv}\em Let be the function $f_{\alpha}(t)$ as defined in (\ref{fa}),  
$x(t)=e^{t}$ and $\Gamma(\alpha)$ as defined in (\ref{sei}); then we have:
\begin{equation}
f_{\alpha}*x = \Gamma(\alpha)x, \quad x(t)=e^t \label{gammaconveq}
\end{equation}
\begin{proof} Function $x(t)\neq 0$ for $t<0$, so we use formula (\ref{aconveq}) for convolution:
\begin{eqnarray*}
(f_{\alpha}*x)(t) & = & \int_{-\infty}^{\infty}f_{\alpha}(\tau)x(t-\tau)d\tau, \quad f_{\alpha}(\tau)=0, \text{ for } \tau<0 \\
         & = & \int_{0}^{\infty}f_{\alpha}(\tau)x(t-\tau)d\tau, \quad x(t-\tau) = e^{t-\tau}=e^{t}e^{-\tau}\\
         & = & \int_{0}^{\infty}f_{\alpha}(\tau)e^{t}e^{-\tau}d\tau\\
         & = & \left(\int_{0}^{\infty}f_{\alpha}(\tau)e^{-\tau}d\tau\right)e^{t}\\
         & = & \left(\int_{0}^{\infty}\tau^{\alpha-1}e^{-\tau}d\tau\right)e^{t}\\
         & = & \Gamma(\alpha)x(t).                
\end{eqnarray*}

\end{proof}

\end{proposition}

\begin{remark}\em Since $f_1$ is the unit step function:
\[
(f_1*x)(t) = \int_{-\infty}^tx(\tau)d\tau = \int_{-\infty}^te^{\tau}d\tau = e^t,
\]
we then have $f_1*x=x$, which implies $\Gamma(1)=1$. 
\end{remark}

\begin{corollary}\em 
Let be $\alpha>0\in\mathbb{R}$, then:
\begin{equation}\label{proprec}
\Gamma(\alpha+1) = \alpha\Gamma(\alpha) \iff \Gamma(\alpha) = \frac{\Gamma(\alpha+1)}{\alpha}
\end{equation}

\begin{proof}

By (\ref{gammaconveq}), we have:

\begin{eqnarray*}
\Gamma(\alpha+1)x & = & f_{\alpha+1}*x, \quad f_{\alpha+1} = \alpha f_{\alpha}*f_1 \\
                & = & \alpha f_{\alpha}*(f_1*x), \quad f_1*x=x \\
                & = & \alpha f_{\alpha}*x \\
                & = & \alpha \Gamma(\alpha)x \implies \\
\Gamma(\alpha+1)\cancel{x} & = &  \alpha \Gamma(\alpha)\cancel{x}                  
\end{eqnarray*}
\end{proof}

\end{corollary}

\begin{remark}\em As explained in \cite{artin}, functional relationship (\ref{proprec}) allow us (recursively) 
to extend the definition of $\Gamma(\alpha)$ for all $\alpha\in\mathbb{R}$, since the integral (\ref{sei}) for 
 $\Gamma(\alpha+1)$ is convergent for $\alpha>-1$.
\end{remark}

We now use (\ref{fabe}) and (\ref{gammaconveq}) to obtain the well-known relationship between $\Gamma$ and $\mathcal{B}$

\begin{corollary}
\begin{equation}\label{sumprop}
\Gamma(\alpha_1)\Gamma(\alpha_2)=\mathcal{B}(\alpha_1,\alpha_2)\Gamma(\alpha_1+\alpha_2) \iff 
\mathcal{B}(\alpha_1,\alpha_2) = \frac{\Gamma(\alpha_1)\Gamma(\alpha_2)}{\Gamma(\alpha_1+\alpha_2)}
\end{equation}

\begin{proof}
\begin{eqnarray*}
\Gamma(\alpha_1)\Gamma(\alpha_2)x & = & \Gamma(\alpha_1)(f_{\alpha_2}*x) \\
                                    & = & f_{\alpha_2}*(\Gamma(\alpha_1)x) \\
                                    & = & f_{\alpha_2}*(f_{\alpha_1}*x)\\
                                    & = & (f_{\alpha_1}*f_{\alpha_2})*x \\
                                    & = & \mathcal{B}(\alpha_1,\alpha_2)(f_{\alpha_1+\alpha_2}*x) \implies \\
\Gamma(\alpha_1)\Gamma(\alpha_2)\cancel{x} & = & \mathcal{B}(\alpha_1,\alpha_2)\Gamma(\alpha_1+\alpha_2)\cancel{x}
\end{eqnarray*}

\end{proof}
\end{corollary}

In the the following Proposition we explore the relationship (\ref{sumprop}) to obtain a generalization of
Proposition~\ref{fab}: 
\begin{proposition}\label{nconvth}\em 
Consider the functions $f_{\alpha_1}, f_{\alpha_2}, \ldots, f_{\alpha_n}$, $\alpha_i>0$,
then:
\begin{equation}\label{nconv}
f_{\alpha_1}*f_{\alpha_2}*\cdots*f_{\alpha_n} = 
\frac{\Gamma(\alpha_1)\Gamma(\alpha_2)\cdots \Gamma(\alpha_n)}
{\Gamma(\alpha_1+\alpha_2+\cdots+\alpha_n)}f_{\alpha_1+\alpha_2+\cdots+\alpha_n}
\end{equation}

\begin{proof}
By induction on $n$. Trivially valid for $n=1$. Suppose valid for $n$, then:
\begin{eqnarray*}
f_{\alpha_1}*f_{\alpha_2}*\cdots*f_{\alpha_n}*f_{\alpha_{n+1}} & = & 
(f_{\alpha_1}*f_{\alpha_2}*\cdots*f_{\alpha_n})*f_{\alpha_{n+1}} \\
 & = & \frac{\prod_{j=1}^n\Gamma(\alpha_j)}
{\Gamma(\sum_{j=1}^n\alpha_j)}\left(f_{\sum_{j=1}^{n}\alpha_j}*f_{\alpha_{n+1}}\right) \\
 & = & \frac{\prod_{j=1}^n\Gamma(\alpha_j)}
{\Gamma(\sum_{j=1}^n\alpha_j)}\mathcal{B}\left(\sum_{j=1}^n\alpha_j, \alpha_{n+1}\right)
f_{\sum_{j=1}^{n+1}\alpha_j}\\
& = & \frac{\prod_{j=1}^n\Gamma(\alpha_j)}
{\cancel{\Gamma(\sum_{j=1}^n\alpha_j)}}\frac{\cancel{\Gamma(\sum_{j=1}^n\alpha_j)}\Gamma(\alpha_{n+1})}
{\Gamma(\sum_{j=1}^{n+1}\alpha_j)}f_{\sum_{j=1}^{n+1}\alpha_j}
\end{eqnarray*}

\end{proof}
\end{proposition}

\begin{corollary}\em Let be $f_{\alpha}^{*n}$ the $n$-fold convolution of $f_{\alpha}$, that is, 
$f_{\alpha}^{*n} = \underbrace{f_{\alpha}*f_{\alpha}*\cdots*f_{\alpha}}_{n\text{ times}}$, $\alpha>0$, then:
\begin{equation}\label{nfold}
f_{\alpha}^{*n} = \frac{\Gamma(\alpha)^n}{\Gamma(n\alpha)}f_{n\alpha}, \quad n\in\mathbb{N}=\{1,2,\ldots \}
\end{equation}

\begin{proof} By assuming $\alpha_1=\alpha_2=\cdots=\alpha_n=\alpha>0$ in (\ref{nconv}).

\end{proof}

Additonally, by taking $\alpha=1/n$ in (\ref{nfold}), we get:
\begin{equation}\label{nfold1}
f_{1/n}^{*n} = \Gamma(1/n)^nf_1, \text{ equivalently, } 
f_{1/n}^{*n}(t) = \begin{cases}0, &t<0 \\ \Gamma(1/n)^n, & t>0 \end{cases}
\end{equation}

For example, since $\Gamma(1/2) = \sqrt{\pi}$, we have ``$f_{1/2}^{*2}(t)=f_{1/2}*f_{1/2}(t)=\pi f_1$" or
\[
f_{1/2}*f_{1/2} = \begin{cases} 0, & t<0 \\ \pi, & t>0 \end{cases}
\]

\end{corollary}

\begin{corollary} \em \label{fa1a}
Let be $0<\alpha<1$, then 
\begin{equation}\label{convext0}
(f_{\alpha}*f_{1-\alpha})(t) = \Gamma(\alpha)\Gamma(1-\alpha)f_1(t)
\end{equation}
Equivalently, from the fact that  ``$\Gamma(\alpha)\Gamma(1-\alpha)= \frac{\pi}{\sin\pi\alpha}$'' we have
 
\begin{equation}\label{convext}
(f_{\alpha}*f_{1-\alpha})(t) = \frac{\pi}{\sin\pi\alpha}f_1(t)
\end{equation}

\begin{proof} Since $\alpha>0$ and $1-\alpha>0$, by using Proposition~\ref{nconvth} we have
\begin{align*}
f_{\alpha}*f_{1-\alpha} & = \frac{\Gamma(\alpha)\Gamma(1-\alpha)}{\Gamma(\alpha+1-\alpha)}f_{\alpha+1-\alpha}\\
                        & = \frac{\Gamma(\alpha)\Gamma(1-\alpha)}{\Gamma(1)}f_{1} \\
                        & = \Gamma(\alpha)\Gamma(1-\alpha)f_{1}
\end{align*}

\end{proof}
\end{corollary}


\begin{corollary}\em\label{hnhncor} Let be $n\ge 2\in\mathbb{N}$ and consider 
\begin{equation}\label{hndef}
h_n = f_{\frac{1}{n}}*f_{\frac{2}{n}}*f_{\frac{3}{n}}*\cdots*f_{\frac{n-1}{n}}
\end{equation}
Then
\begin{equation}\label{hnhneq}
h_n^{*2} = h_n*h_n = \frac{(2\pi)^{n-1}}{n}f_1^{*(n-1)}
\end{equation}

\begin{proof}
\begin{align*}
h_n*h_n & = \left(f_{\frac{1}{n}}*f_{\frac{2}{n}}*\cdots*f_{\frac{n-1}{n}} \right)*
          \left(f_{\frac{1}{n}}*f_{\frac{2}{n}}*\cdots*f_{\frac{n-1}{n}} \right)\nonumber \\
        & = \left(f_{\frac{1}{n}}*f_{\frac{n-1}{n}}\right)*\left(f_{\frac{2}{n}}*f_{\frac{n-2}{n}}\right)* \cdots *
            \left(f_{\frac{n-1}{n}}*f_{\frac{1}{n}}\right) 
            \quad \text{since, by (\ref{convext}), } f_{\frac{k}{n}}*f_{\frac{n-k}{n}} =\frac{\pi}{\sin(\frac{k\pi}{n})}f_1 \\
         & = \left(\frac{\pi}{\sin(\frac{\pi}{n})}f_1\right)*\left(\frac{\pi}{\sin(\frac{2\pi}{n})}f_1\right)*
         \left(\frac{\pi}{\sin(\frac{3\pi}{n})}f_1\right)*\cdots*\left(\frac{\pi}{\sin(\frac{(n-1)\pi}{n})}f_1\right)
          \nonumber \\
          & = \frac{\pi^{n-1}}{\underbrace{\prod_{k=1}^{n-1}\sin\left(\frac{k\pi}{n}\right)}_{=\frac{n}{2^{n-1}}}}f_1^{*(n-1)} 
          \quad \implies\\[0.2cm]
        h_n*h_n & = \frac{(2\pi)^{n-1}}{n}f_1^{*(n-1)}
\end{align*}

\end{proof}
\end{corollary}

The following Propostion is a ``convolutional'' version of the classical Gauss product formula for the Gamma function:

\begin{proposition}\label{convgprodt}\em Let be $\alpha>0\in\mathbb{R}$, $n\ge 2\in\mathbb{N}$, then

\begin{equation}\label{convgprodf}
n^{\alpha-1}f_{\frac{\alpha}{n}}*f_{\frac{\alpha+1}{n}}*f_{\frac{\alpha+2}{n}}*\cdots*f_{\frac{\alpha+n-1}{n}} = 
f_{\frac{1}{n}}*f_{\frac{2}{n}}*\cdots*f_{\frac{n-1}{n}}*f_{\alpha}
\end{equation}

\begin{proof} Define $h_n = f_{\frac{1}{n}}*f_{\frac{2}{n}}*\cdots*f_{\frac{n-1}{n}}$, as in (\ref{hndef}), and 
\begin{equation}\label{gaconv}
g_{\alpha}  = n^{\alpha}h_n*f_{\frac{\alpha}{n}}*f_{\frac{\alpha+1}{m}}*
f_{\frac{\alpha+2}{n}}*\cdots*f_{\frac{\alpha+n-1}{n}}.
\end{equation}
We prove that $g_{\alpha}$ satisfies property (3) of Proposition~\ref{talfa}, that is 
$g_{\alpha+1}=\alpha g_{\alpha}*f_1$. In fact,
\begin{align*}
g_{\alpha+1} & = n^{\alpha+1}h_n*f_{\frac{\alpha+1}{n}}*f_{\frac{\alpha+2}{n}}*\cdots*f_{\frac{\alpha+n-1}{n}}*
\underbrace{f_{\frac{\alpha+n}{n}}}_{=f_{\frac{\alpha}{n}+1}} \\
             & = n^{\alpha+1}h_n*f_{\frac{\alpha+1}{n}}*f_{\frac{\alpha+2}{n}}*\cdots*f_{\frac{\alpha+n-1}{n}}*
f_{\frac{\alpha}{n}+1}, \text{ but } f_{\frac{\alpha}{n}+1} = \frac{\alpha}{n}f_{\frac{n}{\alpha}}*f_1, \text{ so} \\
g_{\alpha+1} & = n^{\alpha+1}h_n*f_{\frac{\alpha+1}{n}}*f_{\frac{\alpha+2}{n}}
             *\cdots*f_{\frac{\alpha+n-1}{n}}*\left(\frac{\alpha}{n}f_{\frac{n}{\alpha}}*f_1\right) \\
             & = \frac{\alpha}{n}n^{\alpha+1}h_n*f_{\frac{n}{\alpha}}*f_{\frac{\alpha+1}{n}}*f_{\frac{\alpha+2}{n}}
             *\cdots*f_{\frac{\alpha+n-1}{n}}*f_1 \\
             & = \alpha \underbrace{n^{\alpha}h_n*f_{\frac{\alpha}{n}}*f_{\frac{\alpha+1}{n}}*f_{\frac{\alpha+2}{n}}
                 *\cdots*f_{\frac{\alpha+n-1}{n}}}_{=g_{\alpha}}*f_1 \implies \\
 g_{\alpha+1} & = \alpha g_{\alpha}*f_1
\end{align*}
On the other hand, by taking $\alpha=1$ in (\ref{gaconv}), we get that $g_1$ is given by 
\begin{align*}
g_1 & =nh_n*\underbrace{f_{\frac{1}{n}}*f_{\frac{2}{n}}*f_{\frac{3}{n}}*\cdots*f_{\frac{n-1}{n}}}_{=h_n}*f_1 \nonumber\\
    & = n(h_n*h_n)*f_1, \quad\text{by (\ref{hnhneq}) } h_n*h_n = \frac{(2\pi)^{n-1}}{n}f_1^{*(n-1)} \implies \\
g_1 & = (2\pi)^{n-1}*f_1^{*n}
\end{align*}

Then we have shown that 
\begin{align*}
g_{\alpha+1} & = \alpha g_{\alpha}*f_1, \text{ and} \\ 
g_1         & = (2\pi)^{n-1}f_1^{*n}.
\end{align*}
So, by Corollary~\ref{talfac}, we can write:
\[
g_{\alpha}=(2\pi)^{n-1}f_1^{*(n-1)}*f_{\alpha}.
\]
Using $g_{\alpha}$ as defined in (\ref{gaconv}) we then get
\begin{align*}
n^{\alpha}h_n*f_{\frac{\alpha}{n}}*f_{\frac{\alpha+1}{n}}*
f_{\frac{\alpha+2}{n}}*\cdots*f_{\frac{\alpha+n-1}{n}} & = (2\pi)^{n-1}f_1^{*(n-1)}*f_{\alpha} \implies\\
n^{\alpha}\underbrace{(h_n*h_n)}_{=\frac{(2\pi)^{n-1}}{n}}*f_1^{*(n-1)}*f_{\frac{\alpha}{n}}*f_{\frac{\alpha+1}{n}}*
f_{\frac{\alpha+2}{n}}*\cdots*f_{\frac{\alpha+n-1}{n}} & = (2\pi)^{n-1}f_1^{*(n-1)}*h_n*f_{\alpha} \implies\\
n^{\alpha-1}\cancel{(2\pi)^{n-1}}*f_1^{*(n-1)}*f_{\frac{\alpha}{n}}*f_{\frac{\alpha+1}{n}}*
f_{\frac{\alpha+2}{n}}*\cdots*f_{\frac{\alpha+n-1}{n}} & = \cancel{(2\pi)^{n-1}}f_1^{*(n-1)}*h_n*f_{\alpha} \implies\\
n^{\alpha-1}f_1^{*(n-1)}*f_{\frac{\alpha}{n}}*f_{\frac{\alpha+1}{n}}*
f_{\frac{\alpha+2}{n}}*\cdots*f_{\frac{\alpha+n-1}{n}} & = f_1^{*(n-1)}*h_n*f_{\alpha} \implies\\
n^{\alpha-1}f_{\frac{\alpha}{n}}*f_{\frac{\alpha+1}{n}}*
f_{\frac{\alpha+2}{n}}*\cdots*f_{\frac{\alpha+n-1}{n}} & = h_n*f_{\alpha} 
\quad\text{(by successive derivation of both sides)}
\end{align*}

\end{proof}

\end{proposition}

\begin{corollary}\em Equation~(\ref{convgprodf}) in Proposition~\ref{convgprodt} implies the so-called Gauss 
Product Formula for the Gamma function:
\begin{equation}\label{gmf}
\Gamma\left(\frac{\alpha}{n}\right)\Gamma\left(\frac{\alpha+1}{n}\right)\cdots\Gamma\left(\frac{\alpha+n-1}{n}\right) = 
\frac{(2\pi)^{\frac{n-1}{2}}}{n^{\alpha-\frac{1}{2}}}\Gamma(\alpha).
\end{equation}
\begin{proof}
Let be $x(t)=e^t$ and consider the convolution of both sides of (\ref{convgprodf}) with $x(t)$:
\[
n^{\alpha-1}f_{\frac{\alpha}{n}}*f_{\frac{\alpha+1}{n}}*f_{\frac{\alpha+2}{n}}*\cdots*f_{\frac{\alpha+n-1}{n}}*x = 
f_{\frac{1}{n}}*f_{\frac{2}{n}}*\cdots*f_{\frac{n-1}{n}}*f_{\alpha}*x
\]
By using the result of Proposition~\ref{gammaconv} we get
\[
n^{\alpha-1}\Gamma\left(\frac{\alpha}{n}\right)\Gamma\left(\frac{\alpha+1}{n}\right)\cdots
\Gamma\left(\frac{\alpha+n-1}{n}\right)\cancel{x}=
\underbrace{\Gamma\left(\frac{1}{n}\right)\Gamma\left(\frac{2}{n}\right)\cdots\Gamma\left(\frac{n-1}{n}
\right)}_{=\left(\frac{(2\pi)^{n-1}}{n}\right)^{1/2}}
\Gamma(\alpha)\cancel{x}
\]

\begin{align*}
n^{\alpha-1}\Gamma\left(\frac{\alpha}{n}\right)\Gamma\left(\frac{\alpha+1}{n}\right)\cdots
\Gamma\left(\frac{\alpha+n-1}{n}\right)& = \left(\frac{(2\pi)^{n-1}}{n}\right)^{1/2}\Gamma(\alpha) \implies\\
\Gamma\left(\frac{\alpha}{n}\right)\Gamma\left(\frac{\alpha+1}{n}\right)\cdots
\Gamma\left(\frac{\alpha+n-1}{n}\right) & = \frac{1}{n^{\alpha-1}}\left(\frac{(2\pi)^{n-1}}{n}\right)^{1/2}\Gamma(\alpha)
\end{align*}
from which Equation~(\ref{gmf}) follows easily.
\end{proof}

\end{corollary}

By the fact that $\Gamma(\alpha)$ is defined for all $\alpha\in\mathbb{R}$, we use 
Proposition~\ref{nconvth} for obtain a result for convolution involving $f_{\alpha}$ 
even in cases where the convolution integral may diverge. 
In the Proposition below we explore the fact that Corollary~\ref{fa1a} may be valid for all 
$\alpha\in\mathbb{R}\backslash\mathbb{Z}$:

\begin{proposition}\em\label{convinv}
Let be $\alpha\in\mathbb{R}\backslash\mathbb{Z}$ then, 
\[
(f_{\alpha}*f_{-\alpha})(t) = -\frac{\pi}{\alpha\sin{(\pi\alpha)}}\delta(t)
\]
Equivalently
\[
\int_0^t \tau^{\alpha-1}(t-\tau)^{-\alpha-1}d\tau = -\frac{\pi}{\alpha\sin{(\pi\alpha)}}\delta(t)
\]
 
\begin{proof} Since $\Gamma(\alpha)$ is well defined for all $\alpha\in\mathbb{R}$, we may extend 
Equation~(\ref{convext0}) for all $\alpha\in\mathbb{R}$, that is
\begin{equation}\label{convext1}
(f_{\alpha}*f_{1-\alpha})(t) = \Gamma(\alpha)\Gamma(1-\alpha)f_1(t), \quad \alpha\in\mathbb{R}.
\end{equation}
Now we consider $\alpha<0$ and use property (3) in Proposition~\ref{talfa} to obtain 
``$\dot{f}_{1-\alpha}=-\alpha f_{-\alpha}$.'' 

Differentiating both sides of (\ref{convext1}) we get:
\begin{eqnarray*}
f_{\alpha}*\dot f_{1-\alpha}       & = & \Gamma(\alpha)\Gamma(1-\alpha)\dot f_1, \quad \dot f_1=\delta \implies\\
 f_{\alpha}*[(-\alpha)f_{-\alpha}] & = &  \Gamma(\alpha)\Gamma(1-\alpha)\delta \implies \\
 f_{\alpha}*f_{-\alpha}            & = &  -\frac{\Gamma(\alpha)\Gamma(1-\alpha)}{\alpha}\delta, \quad
 \Gamma(\alpha)\Gamma(1-\alpha) = \displaystyle\frac{\pi}{\sin\pi\alpha}  \implies\\[0.2cm]
 f_{\alpha}*f_{-\alpha}            & = & -\frac{\pi}{\alpha\sin{(\pi\alpha)}}\delta
\end{eqnarray*}

For example, by taking  $\alpha=-1/2$ we get the following result: 
\[
f_{-1/2}*f_{1/2} = -2\pi\delta,
\]
or
\[
\int_0^t\tau^{-3/2}(t-\tau)^{-1/2}d\tau = -2\pi\delta(t).
\]

\end{proof} 

\end{proposition}

\section{Conclusions}
In this note we obtain some results for the convolution of functions in a set by exploring the relationship with
the Eulerian integrals. The results may be used to derive some known properties of the Gamma function. While 
we considered the real case, we think that the results may be extend for the case of complex Gamma function.

\bibliographystyle{plain}

\end{document}